\documentclass[10pt]{article}
\usepackage[utf8]{inputenc}
\usepackage[T1]{fontenc}
\usepackage{amsmath}
\usepackage{amsfonts}
\usepackage{amssymb}
\usepackage[version=4]{mhchem}
\usepackage{stmaryrd}
\usepackage{hyperref}
\hypersetup{colorlinks=true, linkcolor=blue, filecolor=magenta, urlcolor=cyan,}
\usepackage{bbold}
\usepackage{graphicx}
\usepackage[export]{adjustbox}
\graphicspath{ {./images/} }

\title{New Generalizations of Morrie's Law and the Euler Product Formula }

\author{Carlos A. Pérez Aparicio. ${ }^{1}$}
\date{}

\begin{document}
\maketitle
${ }^{1}$COITIRM. Independent Researcher.

Email: \href{mailto:cperapa@gmail.com}{cperapa@gmail.com}

Abstract. In this study, we derived the infinite product representation of the sinc(z) function by expressing it in a trigonometric form. This approach draws parallels with Morrie's Law and Euler's Product formula, and their generalizations. The results presented in this investigation are completely original and are expected to offer new insights into the topic.

Mathematics Subject Classifications: 33-01, 33-02, 33-04, 33-02,40A20

Key Words and Phrases: Morrie's law, Euler Product sinc(z) function,Viète's formula,trigonometric products.

\section*{1. Introduction}
Morrie's Law is a unique trigonometric identity with an interesting backstory. Physicist Richard Feynman named it after learning the identity from a boy named Morrie Jacobs during his childhood. Feynman carried this knowledge with him throughout his life, and as a result, the identity became known as Morrie's Law.

\begin{equation*}
2^{n} \cdot \prod_{k=0}^{n-1} \cos \left(2^{k} \alpha\right)=\frac{\sin \left(2^{n} \alpha\right)}{\sin (\alpha)} \tag{1}
\end{equation*}

The special trigonometric identity has been explored in the works of [1, 2], which have presented geometric and analytical proofs. In this brief note, we establish a new infinite product expression similar to type (1), akin to the Euler sinc(z) product formula [3] or the Viète formula when  $z=\frac{\pi}{2}$ as

\begin{equation*}
\prod_{k=0}^{\infty} \cos \left(2^{-k-1} z\right)=\frac{\sin (z)}{z} \tag{2}
\end{equation*}

\section*{2. New Generalizations of the Morrie's law.}
It is possible to generalize (1) such that $\{m, n\} \gg 0,\{q, n, m\} \in Z$ as follows:

\begin{equation*}
\prod_{k=m}^{n-1} q \sin \left(z q^{k}\right) \csc \left(z q^{k+1}\right)=q^{n-m} \sin \left(z q^{m}\right) \csc \left(z q^{n}\right) \tag{3}
\end{equation*}

When $m=0$ and $q=2$ one gets formula (1)

\begin{equation*}
\prod_{k=0}^{n-1} \cos \left(2^{k} z\right)=2^{-n} \csc (z) \sin \left(2^{n} z\right) \tag{4}
\end{equation*}

\textbf{Proof.}
The proof  (3) is by induction. Let $n \rightarrow n+1, \mathrm{z} \in \mathbb{C} \quad$ so multiplying both sides of the equation (3) by

\begin{equation*}
q \sin \left(z q^{n}\right) \csc \left(z q^{n+1}\right) \tag{5}
\end{equation*}

\begin{align*}
&\prod_{k=m}^{n-1} q \sin \left(z q^{k}\right) \csc \left(z q^{k+1}\right) \left(q \sin \left(z q^{n}\right) \csc \left(z q^{n+1}\right)\right) \\
&= \left(q \sin \left(z q^{n}\right) \csc \left(z q^{n+1}\right)\right) \left(q^{n-m} \sin \left(z q^{m}\right) \csc \left(z q^{n}\right)\right) \tag{6}
\end{align*}

After performing some basic simplifications, we obtain the desired result.

\begin{equation*}
\prod_{k=m}^{n} q \sin \left(z q^{k}\right) \csc \left(z q^{k+1}\right)=q^{-m+n+1} \sin \left(z q^{m}\right) \csc \left(z q^{n+1}\right) \tag{7}
\end{equation*}

\section*{3. Sinc(z) function formulas and New Generalizations of the Morrie's law type.}
We introduce two novel generalizations of the infinite product, resembling Euler product formulas, applicable when $q > 2,\quad q \in \mathbb{Z}$ and , $\mathrm{z} \in \mathbb{C}$. The formulations are presented as follows:

\begin{align*}
&\prod_{k=0}^{\infty} q^{-q^{k}(k(q-1)+q)} z^{(q-1) q^{k}} \cos^{-q^{k+1}}\left(z q^{-k-1}\right) \cos^{q^{k}}\left(z q^{-k}\right) \\
&\quad \times \tan^{-q^{k+1}}\left(z q^{-k-1}\right) \tan^{q^{k}}\left(z q^{-k}\right) = \frac{\sin(z)}{z} \tag{8}
\end{align*}

\begin{align*}
&\prod_{k=0}^{\infty} \frac{\sin \left(z q^{-k}\right) \csc \left(z q^{-k-1}\right)}{q} = \frac{\sin (z)}{z} \tag{9}
\end{align*}
The formula( 9) can also be written as
\begin{equation*}
\prod _{k=0}^{\infty } \sum _{n=1}^q \frac{\cos \left(\frac{(2 n-1) z}{(2 q)^{k+1}}\right)}{q}\tag{10}
\end{equation*}

Equivalent formulas (9) have been derived elsewhere using both Chebyshev polynomials and Fourier Transform, as demonstrated by Nishimura [5] and Kent E. Morrison [4].Additionally formula of (8) for the case when q=2 has been established in [6].Formula (9) can also be generalized as:

\begin{align*}
\prod_{k=0}^{\infty} \prod_{\mathrm{k} 1=1}^{q} \frac{\Gamma\left(\frac{\mathrm{k} 1}{q}\right)^{2}}{\Gamma\left(\frac{\mathrm{k} 1-\frac{q^{-k} z}{\pi}}{q}\right) \Gamma\left(\frac{\frac{z q^{-k}}{\pi}+\mathrm{k} 1}{q}\right)}=\frac{\sin (z)}{z} \tag{11}
\end{align*}

The formula (11) is presented, purely out of curiosity; its proof will be considered in a subsequent work.

Next, we introduce a generalization of the previously mentioned equations (8) and (9), resembling the patterns found in Morrie's law formula, with the conditions $\{m, n\} \gg 0,\{q, n, m\} \in Z, \mathrm{z} \in \mathbb{C}$

\begin{align*}
& \prod_{k=m}^{n-1} q^{q^{-k-1}(k(-q)+k+1)} z^{-\left((q-1) q^{-k-1}\right)} \\
& \quad \times \cos ^{q^{-k}}\left(z q^{k}\right) \tan ^{q^{-k}}\left(z q^{k}\right) \cot ^{q^{-k-1}}\left(z q^{k+1}\right)  \tag{12} \\
& = q^{n q^{-n}-m q^{-m}} z^{q^{-n}-q^{-m}} \\
& \quad \times \cos ^{q^{-m}}\left(z q^{m}\right) \tan ^{q^{-m}}\left(z q^{m}\right) \cos ^{-q^{-n}}\left(z q^{n}\right) \cot ^{q^{-n}}\left(z q^{n}\right)
\end{align*}
\begin{align*}
\prod _{k=m}^{n-1} q \sin \left(z q^k\right) \csc \left(z q^{k+1}\right)=q^{n-m} \sin \left(z q^m\right) \csc \left(z q^n\right) \tag{13}
\end{align*}
\textbf{Proof.}
The proof (12) being by induction. Let $n \rightarrow n+1, \mathrm{z} \in \mathbb{C} \quad$ so multiplying both sides of the equation (11) by

\begin{equation*}
q^{q^{-n-1}(n(-q)+n+1)} z^{-\left((q-1) q^{-n-1}\right)} \cos ^{q^{-n}}\left(z q^{n}\right) \tan ^{q^{-n}}\left(z q^{n}\right) \csc ^{q^{-n-1}}\left(z q^{n+1}\right) \tag{14}
\end{equation*}

Multiplying the right-hand side and simplifying we get the formula (11). We achieve the desired outcome.

\textbf{Example 1.} Infinitely nested radicals.

Using several values for $n$ integers $\frac{2 n \sin \left(\frac{\pi }{2 n}\right)}{\pi }$ in formula (8) get the following table:

\begin{tabular}{|c|c|}
\hline
Product $\infty$ & z\\
\hline
$\frac{\sqrt{2-\sqrt{2}} \pi ^3}{16384 \left(2-\sqrt{\sqrt{\sqrt{2}+2}+2}\right)^2}\text{..}$ & $\frac{4 \sqrt{2-\sqrt{2}}}{\pi }$ \\
\hline
$\frac{\sqrt{\left(1-\frac{2}{\sqrt{5}}\right) \left(\frac{\sqrt{5}}{8}+\frac{5}{8}\right)} \pi ^3}{256000 \left(\frac{1}{4} \left(\sqrt{2}-1\right) \sqrt{\frac{1}{2} \left(\sqrt{2}+2\right) \left(\sqrt{5}+5\right)}-\frac{1}{8} \sqrt{2-\sqrt{2}} \left(\sqrt{2}+1\right) \left(\sqrt{5}-1\right)\right)^4}$\text{..} & $\frac{5 \left(\sqrt{5}-1\right)}{2 \pi }$ \\
\hline
$\frac{\pi ^3}{1024 \sqrt{2} \left(2-\sqrt{\sqrt{2}+2}\right)^2}\text{..}$ & $\frac{2 \sqrt{2}}{\pi}$ \\
\hline
$\frac{\pi ^3}{110592 \left(\frac{1}{4} \sqrt{2-\sqrt{2}} \left(\sqrt{2}+1\right)-\frac{1}{4} \left(\sqrt{2}-1\right) \sqrt{3 \left(\sqrt{2}+2\right)}\right)^4}\text{..}$ & $\frac{3}{\pi }$ \\
\hline
\end{tabular}
\begin{center}
Table 1. 
\end{center}
\textbf{Example 2.} for $q=3$  in equation(8)  $ \mathrm{z} \in \mathbb{C}$
\begin{equation*}
\prod _{k=0}^{\infty } 3^{-3^k (2 k+3)} z^{2\ 3^k} \cos ^{-3^{k+1}}\left(3^{-k-1} z\right) \cos ^{3^k}\left(3^{-k} z\right) \tan ^{-3^{k+1}}\left(3^{-k-1} z\right) \tan ^{3^k}\left(3^{-k} z\right)=\frac{\sin (z)}{z}\tag{15}
\end{equation*}
             
\section*{References}

\begin{enumerate}
    \item Beyer, W. A., Louck, J. D., and Zeilberger, D., A Generalization of a Curiosity that Feynman Remembered All His Life, Math. Mag. 69, 43-44, 1996. (JSTOR).
    
    \item Moreno, S. G., García-Caballero, E. M., "A Geometric Proof of Morrie's Law". American Mathematical Monthly, vol. 122, no. 2 (February 2015), p. 168 (JSTOR).
    
    \item Euler, Leonhard (1738). "De variis modis circuli quadraturam numeris proxime exprimendi".
    
    \item Morrison, K. E., Cosine Products, Fourier Transforms, and Random Sums, The American Mathematical Monthly, 102:8, 716-724, DOI: \url{10.1080/00029890.1995.12004647}.
    
    \item Nishimura, R., A generalization of Viète's infinite product and new mean iterations. The Australian Journal of Mathematical Analysis and Applications, 13(1), 2016.
    
    \item Milgram, M., A Curious Trigonometric Infinite Product in Context. arXiv:2303.08628 [\url{http://math.GM}].
\end{enumerate}

\end{document}